\documentclass[12pt]{amsart}
\usepackage{amssymb,latexsym, amsmath, amsxtra}
\usepackage{graphicx}
\begin{document}

\title{Moments of zeta and correlations of divisor-sums: II}
 
\author{Brian Conrey}
\address{American Institute of Mathematics, 360 Portage Ave, Palo Alto, CA 94306, USA and School of Mathematics, University of Bristol, Bristol BS8 1TW, UK}
\email{conrey@aimath.org}
\author{Jonathan P. Keating}
\address{School of Mathematics, University of Bristol, Bristol BS8 1TW, UK}
\email{j.p.keating@bristol.ac.uk}

\thanks{We gratefully acknowledge support under EPSRC Programme Grant EP/K034383/1
LMF: L-Functions and Modular Forms.  Research of the first author was also supported by the American Institute of Mathematics and by a grant from the National Science Foundation. JPK is grateful for the following additional support: a grant from the Leverhulme Trust, a Royal Society Wolfson Research 
Merit Award, a Royal Society Leverhulme Senior Research Fellowship, and a grant from the Air Force Office of Scientific Research, Air Force Material Command, 
USAF (number FA8655-10-1-3088). He is also pleased to thank the American Institute of Mathematics for hospitality during a visit where this work started.}

\date{\today}

\begin{abstract}
This is part II of our
examination  of the second  and fourth moments and shifted moments of the Riemann zeta-function on the critical line using long Dirichlet 
polynomials and divisor correlations. 

\end{abstract}

\maketitle

\section{Introduction}
In part I, see [CK],  we  
completed the analysis of the second moment of the Riemann zeta-function using the long Dirichlet polynomial method of Goldston and Gonek
[GG] and we initiated the study of the fourth moment by this approach. 
In particular we calculated the contributions from the off-diagonal terms arising from coefficient correlations
of the form $\sum_{n\le X} d(n)d(n+h)$ and  
identified the terms that are missed in this approach.  In this paper we show how
to evaluate these new terms that were missing and in doing so we introduce a new technique that is a discrete analog of the circle method. 
This analysis gives a concrete introduction to how we will approach higher moments through this circle method approach. In a subsequent paper 
we will show how to obtain the ``full-moment'' conjecture for the $2k$th moment of $\zeta(s)$ on the critical line, i.e. the full polynomial of degree
$k^2$ which comprises the main term.  
The idea for this method originates in the work of Bogomolny and Keating; see [BK].

Thus, we will calculate the contribution of what we call the type II sums (after [BK]) which arise in the evaluation of 
\begin{eqnarray*}
\int_0^\infty \psi\left(\frac t T\right) \sum_{m\le X}\frac{\tau_{\alpha,\beta}(m)}{m^{s}}
\sum_{m\le X}\frac{\tau_{\gamma,\delta}(n)}{n^{1-s}}~dt 
\end{eqnarray*}
where $s=1/2+it$ and $\tau_{\alpha,\beta}(n)=\sum_{de=n}d^{-\alpha}e^{-\beta}$. (See [CK] for further
notation and introduction.) To describe the type II sums we observe that 
integrating term-by-term we find that the above is
$$T\sum_{m,n\le X} \frac{\tau_{\alpha,\beta}(m)\tau_{\gamma,\delta}(n)}{\sqrt{mn}} \hat \psi(\frac{T}{2\pi}\log(m/n))
=T\mathcal D +T\mathcal O +\mathcal E$$
where $\mathcal D $ is the diagonal
$$\mathcal D =\hat{\psi}(0)  \sum_{n\le X} \frac{\tau_{\alpha,\beta}(n)\tau_{\gamma,\delta}(n)}{n};$$
$\mathcal O $ is the off-diagonal 
$$\mathcal O=\sum_{m\ne n\atop
0<|m-n|<m/\tau} \frac{\tau_{\alpha,\beta}(m)\tau_{\gamma,\delta}(n)}{\sqrt{mn}} \hat \psi(\frac{T}{2\pi}\log(m/n));$$
and $\mathcal E \ll T^{\epsilon}$ is an error term; here $\tau=T^{1-\epsilon}$
and the Fourier transform is defined by
$$\hat{\psi}(v)=\int_{\mathbb R} \psi(u) e(uv)~du$$
where $e(x)=\exp(2\pi i x)$.

If we evaluate  $\mathcal O$ here in the traditional manner, eg. as in [GG], we would now solve the 
shifted convolution problem which consists of  evaluating
$$\sum_{n\le x} \tau_{\alpha,\beta}(n)\tau_{\gamma,\delta}(n+h)$$  
and summing by parts. This analysis was carried out in I. Here we use  a new approach.
We first make use of the fact that $\tau_{\alpha,\beta}$ and $\tau_{\gamma,\delta}$ are convolutions to write 
$$\mathcal O=\sum_{m_1m_2,n_1n_2 \le X\atop
0<|m_1m_2-n_1n_2|<m_1m_2/\tau} \frac{m_1^{-\alpha}m_2^{-\beta} n_1^{-\gamma}n_2^{-\delta}}{m_1m_2} \hat \psi(\frac{T}{2\pi}\log((n_1n_2)/(m_1m_2)).$$
Now we embark on a discrete analog of the circle method which basically consists of approximating a ratio, say $m_1/n_1$ by 
a rational number with a small denominator, say $M/N$, and then sum all of the terms with $m_1/n_1$ close to $M/N$.

 To this end  
we introduce a parameter $Q$ and subdivide the interval $[0,1]$ into Farey intervals associated with  the fractions $M/N$ with 
$1\le M\le N\le Q$ and $(M,N)=1$ from the Farey sequence $\mathcal F_Q$. The Farey interval $\mathcal M_{M,N}$ determined by the fraction $M/N$ is defined to be 
$$\mathcal M_{M,N} = \left[ \frac MN -\frac {M+M''}{N+N''}, \frac MN +\frac {M+M'}{N+N'}\right)$$
where 
$\frac{M''}{N''}, \frac M N , \frac {M'}{N'}$ 
are three consecutive terms in the Farey sequence $\mathcal F_Q$. Now given such an $M$ and $N$ we 
sum over the terms  $m_1$ and $n_1$ for which $m_1/n_1\in \mathcal{M}_{M,N}$; for such a pair we 
  define 
$$h_1 := m_1 N - n_1 M.$$
The possible range of $h_1$ 
may be computed by 
$$ |h_1|  = \left| \frac {m_1}{n_1}-\frac  M N \right| n_1 N \le  \left(\frac MN -\frac {M+M''}{N+N''}\right)n_1 N
=\frac{n_1}{N+N''} \approx \frac{n_1}{Q}$$
since adjacent denominators satisfy $Q<N+N''<2Q$. In general, the rapid decay of $\hat{\psi}$ governs the 
range of $h_1$ and $h_2$ defined below. 

Also, we note that if $Q$ is not too large then $m_1/n_1\in \mathcal M_{M,N}$ implies that $n_2/m_2\in \mathcal M_{M,N}$
as well. This is because the distance from $m_1/n_1$ to $n_2/m_2$ is 
$$ \left| \frac {m_1}{n_1} -\frac{n_2}{m_2}\right|=\frac{|m_1m_2-n_1n_2|}{ n_1m_2}
\le \frac{m_1m_2}{\tau n_1 m_2}\le \frac{1}{\tau}.$$
On the other hand 
$$\left| \frac M N -\frac {M'}{N'}\right| \gg \frac{1}{Q^2}$$
so if $Q^2 =o( \tau)$ then our assertion follows.

Now we define $$h_2: = m_2 M -n_2 N.$$
We have
$$ m_1m_2 MN -n_1n_2 MN = h_1m_2M+h_2m_1 N -h_1 h_2$$
so that 
$$
\frac{m_1 m_2 -n_1 n_2}{m_1m_2}=\frac{h_1}{m_1N}+\frac{h_2}{m_2 M} -\frac{h_1h_2}{m_1m_2MN}
$$
and 
$$\log\frac{n_1n_2}{m_1m_2}=\frac{h_1}{m_1N}+\frac{h_2}{m_2 M} +O\big(\frac{h_1h_2}{m_1m_2MN}\big).
$$
The error term is negligible so we 
have now  arranged the sum as 
$$\sum_{M\le N\le Q\atop (M,N)=1}
\sum_{h_1,h_2 }
\sum_{{m_1m_2\le X }\atop {(*_1), (*_2) }}
\frac{m_1^{-\alpha}m_2^{-\beta}n_1^{-\gamma}n_2^{-\delta}}{m_1m_2}
\hat \psi \left(\frac{Th_1}{2\pi m_1N}+\frac{Th_2}{2\pi m_2 M} \right)
$$
 where 
$$ (*_1): m_1N-n_1M=h_1 \qquad \qquad \mbox{ and } \qquad \qquad  (*_2): m_2M-n_2N=h_2
$$
Note that for a given $m_1, n_1$ and $h_1$ the condition $(*_1)$ implies that $m_1/n_1 \in \mathcal M_{M,N}$
so we don't need to write that condition. 

\section{Smoothing the sums over $M$ and $N$ }
We introduce another smooth weight function $\phi(y)$, which is an approximation to the characteristic function 
$\chi_{(0,1]}(y)$ to help with the summation over $M$ and $N$.  
We will then have sums of the form  
$$S_Q(\xi,\eta):=\sum_{1\le M\le N 
\atop (M,N)=1 } \phi\left( \frac M Q\right)  \phi\left( \frac N Q\right) M^{-1-\xi} N^{-1-\eta}$$
 for a finite set of choices of $\xi$ and $\eta$ which are of the form 
$$\epsilon_1\alpha+\epsilon_2\beta
+\epsilon_3\gamma+\epsilon_4\delta$$ 
where the $\epsilon_i \in \{\-1,0,1\}$. 
We require that 
$$\phi(y) =\frac{1}{2\pi i} \int_{(1)} \tilde{\phi} (s) y^{-s} ~ds$$
where $\tilde{\phi}(s)$ has the properties that
$$\operatornamewithlimits{Res}_{s=0}\tilde{\phi}(s)=1 \qquad \mbox{ and } \qquad \tilde{\phi}(\xi)=0$$
for all of the eligible values of $\xi$ that arise, and that $\tilde{\phi}(s)$ is analytic in $\Re s\ge -1/2$ 
and has rapid decay vertically in this region.  
In practice $S_Q(\xi,\eta)$ will be combined   with
$S_Q(\eta,\xi)$ to obtain
\begin{eqnarray*}
S_Q(\xi,\eta)+S_Q(\eta, \xi)&=& \phi\left(\frac 1 Q\right)^2 +\sum_{(M,N)=1}\phi\left( \frac M Q\right) 
 \phi\left( \frac N Q\right) M^{-1-\xi} N^{-1-\eta}
\end{eqnarray*}
The   second term is 
\begin{eqnarray*}&&
\sum_{d }\frac{\mu(d)}{d^{2+\xi+\eta}} \sum_{M } \phi\left(\frac{Md}{Q}\right) M^{-1-\xi}
\sum_{N } \phi\left(\frac{Nd}{Q}\right) M^{-1-\eta}
 \\
&=& \sum_{d }\frac{\mu(d)}{d^{2+\xi+\eta}}
\bigg(\frac{1}{2\pi i}\int_{(1)} \tilde{\phi}(w) \zeta(w+1+\xi)  \left(\frac{Q}{d}\right)^{w} ~dw\bigg)
\bigg(\frac{1}{2\pi i}\int_{(1)} \tilde{\phi}(z) \zeta(z+1+\eta)  \left(\frac{Q}{d}\right)^{z} ~dz\bigg) .
\end{eqnarray*}
The first integral is $= \zeta(1+\xi)+O((Q/d)^{-1/3})$ as can be seen by moving the path of integration to the left to $\Re w=-1/3$
and accounting for the residue at the pole $w=0$; note that since $\tilde{\phi}(-\xi)=0$, there is no pole at $w=-\xi$.
Thus, altogether we have 
\begin{eqnarray} \label{eqn:Ssum}
S_Q(\xi,\eta)+S_Q(\eta, \xi)&=& \phi\left(\frac 1 Q\right)^2 + \frac{\zeta(1+\xi)\zeta(1+\eta)}{\zeta(2+\xi+\eta)}  +O(Q^{-1/3}).
\end{eqnarray}

\section{The case of $h_2=0$} 
We remark first of all that the terms with $h_1=h_2=0$ are precisely the diagonal terms. Now we consider what happens if $h_2=0$ and $h_1\ne0$.
We call this a ``semi-diagonal'' term after [BK].

If $h_2=0$ then $m_2M=n_2N$. Since $(M,N)=1$ it follows that
$m_2=N\ell$ and $n_2=M\ell$ for some $\ell$.
Thus we have
$$\sum_{M\le N\atop (M,N)=1}\phi\left( \frac M Q\right)  \phi\left( \frac N Q\right)
\sum_{ h_1}
\sum_{{m_1, n_1 ,\ell }\atop {(*_1) \atop n_1\ge |h_1|Q }}
\frac{m_1^{-\alpha}(N\ell)^{-\beta}n_1^{-\gamma}(M\ell)^{-\delta}}{m_1m_2}
\hat \psi \left(\frac{Th_1}{2\pi m_1N} \right)
$$
where
$$(*_1): m_1N-n_1M=h_1 .$$
We replace $m_1$ by a smooth variable $u_1$ and $n_1$ by $m_1N/M$.
  We have $u_1\ell N=m_1m_2\le X$ and so our sum is 
$$\sum_{M\le N\atop (M,N)=1} \phi\left( \frac M Q\right)  \phi\left( \frac N Q\right) M^{-\delta+\gamma-1}N^{-\beta-\gamma-1}
\sum_{h_1}
\sum_{ \ell  }\ell^{-1-\beta-\delta} \int_{u_1\ell \le \frac{X}{ N} } u_1^{-1-\alpha-\gamma}
\hat \psi \left(\frac{Th_1}{2\pi u_1N} \right)~du_1
$$
We save the term with $h_1=0$ for later and we group the terms with $h_1$ and $-h_1$ together and use $\hat \psi(-v)=\overline {\hat \psi(v)}$.  
We make the substitution $v_1=\frac{Th_1}{2\pi u_1N}$ in the integral and switch the integral over $v_1$ with the sum over $h_1$ and $\ell$. 
Then (with $h_1>0$) we have that  
$$\frac{\ell NTh_1}{2\pi v_1N}=u_1\ell N\le X
$$
implies that $$\ell h_1 \le \frac{2\pi Xv_1}{T}.$$
Thus we have
$$\sum_{M\le N\atop (M,N)=1} \phi\left( \frac M Q\right)  \phi\left( \frac N Q\right)  M^{-\delta+\gamma-1}N^{\alpha-\beta-1}
 \int_{0}^\infty v_1^{-1+\alpha+\gamma} (2\Re \hat{\psi}(v_1)) \sum_{h_1\ell \le \frac{2\pi Xv_1}{T}}
h_1^{-\alpha-\gamma} \ell^{-1-\beta-\delta} ~dv_1.
$$
The sum over $h_1$ and $\ell$   is 
$$\frac{1}{2\pi i} \int_{(2)} \zeta(s+1+\beta+\delta )\zeta(s+\alpha+\gamma)  \left(\frac {2\pi v_1X}{T}\right)^s \frac{ds}{s} 
$$  
Together with the integral over $v_1$ this is 
\begin{eqnarray*}
\int_{0 }^{\infty} v_1^{-1+\alpha+\gamma}
\hat \psi(v_1)
\frac{2}{2\pi i} \int_{(2)} \zeta(s+1+\beta+\delta )\zeta(s+\alpha+\gamma)  \left(\frac {2\pi v_1X}{T}\right)^s \frac{ds}{s} ~dv_1.
\end{eqnarray*}
Now, as we've seen before, if $\Re s>0$ then 
\begin{eqnarray*}
\int_0^\infty (2\Re \hat {\psi}(v)) v^s~dv=\chi(1-s) \int_0^\infty \psi(t)t^{-s}~dt.
\end{eqnarray*}
Thus, the above is 
 \begin{eqnarray*}
\int_{0 }^{\infty} t^{-1-\alpha-\gamma}
 \psi(t)
\frac{2}{2\pi i} \int_{(2)} \zeta(s+1+\beta+\delta )\zeta(1-s-\alpha-\gamma)  \left(\frac {2\pi X}{tT}\right)^s \frac{ds}{s} ~dt
\end{eqnarray*}
We move the $s$-path left to $\Re s=-1/2$, thus crossing the poles at $s=0$, $s=-\alpha-\gamma$ and $s=-\beta-\delta$. Thus the 
above is 
\begin{eqnarray*}&&
\int_{0 }^{\infty} t^{-1-\alpha-\gamma}
 \psi(t)\bigg(\zeta(1+\beta+\delta )\zeta(1-\alpha-\gamma) -\frac{\zeta(1-\alpha+\beta-\gamma+\delta) \left(\frac {2\pi X}{tT}\right)^{-\beta-\delta}}
{\beta+\delta}\\
&& \qquad \qquad \qquad \qquad \qquad \qquad 
+ \frac{\zeta(1-\alpha+\beta-\gamma+\delta) \left(\frac {2\pi X}{tT}\right)^{-\alpha-\gamma}}
{\alpha+\gamma}\bigg)~dt
\end{eqnarray*}
and altogether we have
\begin{eqnarray*}&&
\sum_{M\le N\atop (M,N)=1} \phi\left( \frac M Q\right)  \phi\left( \frac N Q\right) M^{-\delta+\gamma-1}N^{\alpha-\beta-1}
\\&&\qquad \times
\int_{0 }^{\infty} t^{-1-\alpha-\gamma}
 \psi(t)\bigg(\zeta(1+\beta+\delta )\zeta(1-\alpha-\gamma) -\frac{\zeta(1-\alpha+\beta-\gamma+\delta) \left(\frac {2\pi X}{tT}\right)^{-\beta-\delta}}
{\beta+\delta}\\
&& \qquad \qquad \qquad \qquad \qquad \qquad 
+ \frac{\zeta(1-\alpha+\beta-\gamma+\delta) \left(\frac {2\pi X}{tT}\right)^{-\alpha-\gamma}}
{\alpha+\gamma}\bigg)~dt
\end{eqnarray*}

All of the above is predicated on $m_1/n_1<1$. The contribution from the terms where $n_1<m_1$ will
be exactly as above but with the quadruple $(\alpha,\beta,\gamma,\delta)$ replaced with $(\gamma,\delta,\alpha,\beta)$.
In particular, $\alpha+\gamma$ will be replaced by $\beta+\gamma$ prior to summing over $M$ and $N$. 
This will give another term 
\begin{eqnarray*}&&
\sum_{M\le N\atop (M,N)=1} \phi\left( \frac M Q\right)  \phi\left( \frac N Q\right) M^{\alpha-\beta-1}N^{-\delta+\gamma-1}
\\&&\qquad \times
\int_{0 }^{\infty} t^{-1-\beta-\delta}
 \psi(t)\bigg(\zeta(1+\beta+\delta )\zeta(1-\alpha-\gamma) -\frac{\zeta(1-\alpha+\beta-\gamma+\delta) \left(\frac {2\pi X}{tT}\right)^{-\beta-\delta}}
{\beta+\delta}\\
&& \qquad \qquad \qquad \qquad \qquad \qquad 
+ \frac{\zeta(1-\alpha+\beta-\gamma+\delta) \left(\frac {2\pi X}{tT}\right)^{-\alpha-\gamma}}
{\alpha+\gamma}\bigg)~dt
\end{eqnarray*}

Now we consider what happens when $h_1=0$ and $h_2\ne 0$. These terms will contribute the ``complements''
to the above two expressions so that we will be in 
the situation described in  (\ref{eqn:Ssum}) 
  and so we can execute the sums over $M$ and $N$ as described there,
replacing the sums over $M$ and $N$ by ratios of zeta functions with small error terms.  
Thus, we obtain 
\begin{eqnarray*}&&
\int_{0 }^{\infty} t^{-1-\alpha-\gamma}
 \psi(t)\bigg(\frac{\zeta(1+\beta+\delta )\zeta(1-\alpha-\gamma)\zeta(1-\gamma+\delta)\zeta(1-\alpha+\beta)}{
\zeta(2-\alpha +\beta-\gamma +\delta)}
\\
&& \qquad \qquad \qquad \qquad  
 -\left(\frac {2\pi X}{tT}\right)^{-\beta-\delta}\frac{\zeta(1-\alpha+\beta-\gamma+\delta)\zeta(1-\alpha+\beta)\zeta(1-\gamma+\delta)}{
(\beta+\delta) ~ \zeta(2-\alpha+\beta-\gamma+\delta)}\\
&& \qquad \qquad \qquad \qquad \qquad \qquad 
+  \left(\frac {2\pi X}{tT}\right)^{-\alpha-\gamma}\frac{\zeta(1-\alpha+\beta-\gamma+\delta)\zeta(1-\alpha+\beta)\zeta(1-\gamma+\delta)}{
(\alpha+\gamma) ~ \zeta(2-\alpha+\beta-\gamma+\delta)}\bigg)~dt
\end{eqnarray*}
and the complimentary term with $\alpha+\gamma$ replaced by $\beta+\delta$ and vice-versa.

This is identical with one of the  one-swap terms identified by descending as previously described. 

There are further semi-diagonal terms. 
If we do the exact same analysis as throughout this entire section but now focusing on the ratio 
$m_1/n_2$ instead of $m_1/n_1$ then the effect will be to switch the roles of $\gamma$ and $\delta$ in the expression above.
Then we end up with two more terms and a total of four terms. These terms are identical with the 
four terms obtained by the ``descent'' method described in section 8 of [CK].

A question of whether we have over-counted some terms may arise. But the ``duplicate'' terms for which $m_1/n_1 \in 
\mathcal M_{M,N}$ and simultaneously $m_1/n_2 \in \mathcal M_{M',N'}$ with $N\le Q$ and $N'\le Q$ 
  contribute an insignificant amount to the total and so may be regarded as part of the error term. 

\section{The case of $h_1h_2\ne 0$}
Now we
consider 
$$\sum_{M\le N \atop (M,N)=1}\phi\left(\frac{M}{Q}\right)\phi\left(\frac{N}{Q}\right)
\sum_{h_1h_2\ne 0 }
\sum_{{m_1m_2\le X }\atop {(*_1), (*_2) }}
\frac{m_1^{-\alpha}m_2^{-\beta}n_1^{-\gamma}n_2^{-\delta}}{m_1m_2}
\hat \psi \left(\frac{Th_1}{2\pi m_1N}+\frac{Th_2}{2\pi m_2 M} \right).
$$
In this case we  have a  bound for $h_2$ similar to that for $h_1$:
$$|h_2| \ll \frac {m_2}{Q} \ll \frac{n_2 M}{QN}.$$
In particular, we have 
$$|h_1 h_2| \ll \frac{n_1n_2M}{Q^2N} \ll \frac{X}{Q^2}.
$$
 
Now we replace the sums over $m_1,m_2,n_1,n_2$ subject to $(*_1)$ and $(*_2)$ by their averages. 
As before, we replace $m_1$ by $u_1$ and now we replace $m_2$ by $u_2$. We replace $n_1$ and $n_2$ by 
$u_1N/M$ and $u_2M/N$ respectively.  
  We then have
$$\sum_{M\le N \atop (M,N)=1}\phi\left(\frac{M}{Q}\right)\phi\left(\frac{N}{Q}\right)M^{\gamma-\delta-1}N^{\delta-\gamma-1}
\sum_{h_1h_2\ne 0  }
\int_{u_1u_2\le X  }
u_1^{-\alpha-\gamma-1}u_2^{-\beta-\delta-1}
\hat \psi \left(\frac{Th_1}{2\pi u_1N}+\frac{Th_2}{2\pi u_2 M} \right) ~du_1~du_2.
$$ 
Now there are four cases to consider according to the four sign choices of $h_1$ and $h_2$.
We make the substitutions  
$$v_1= \frac{T|h_1|}{2\pi u_1N} \qquad  \mbox{ and } \qquad   v_2=\frac{T|h_2|}{2\pi u_2 M}$$
and move the sums over $h_1$ and $h_2$ to the inside. 
The condition $u_1u_2\le X$ implies that 
$$ \frac{T^2|h_1h_2|}{4\pi^2 v_1v_2MN}  =u_1u_2 \le X
 $$
or
$$ |h_1h_2| \le \frac{4\pi^2 XMNv_1v_2}{T^2} .
$$
We get 
\begin{eqnarray*} &&
\big(\frac{T}{2\pi}\big)^{-\alpha-\beta-\gamma-\delta}\sum_{M\le N \atop (M,N)=1}\phi\left(\frac{M}{Q}\right)\phi\left(\frac{N}{Q}\right)
M^{\gamma+\beta-1}N^{\delta+\alpha-1}\iint_{v_1,v_2}
v_1^{\alpha+\gamma-1}v_2^{\beta+\delta-1}
\sum_{0<|h_1h_2| \le \frac {4\pi^2 XMNv_1v_2}{T^2}}
\\
&&\qquad \times 
h_1^{-\alpha-\gamma} h_2^{-\beta-\delta}
\left(\hat \psi(v_1+v_2)+\hat \psi(v_1-v_2)+\hat \psi(-v_1+v_2)+\hat \psi(-v_1-v_2)\right) ~dv_1~dv_2.
\end{eqnarray*}
Using
$$\hat\psi(v_1+v_2)=\int_0^\infty \psi(t)e(t(v_1+v_2))~dt$$
we see that 
$$\hat \psi(v_1+v_2)+\hat \psi(v_1-v_2)+\hat \psi(-v_1+v_2)+\hat \psi(-v_1-v_2)
=\int_0^\infty \psi(t)\big(e(tv_1)+e(-tv_1)\big)\big(e(tv_2)+e(-tv_2)\big)~dt;
$$
Also
$$\sum_{h_1h_2\le \frac {4\pi^2 XMNv_1v_2}{T^2}}h_1^{-\alpha-\gamma} h_2^{-\beta-\delta}
=\frac{1}{2\pi i}\int_{(2)} \zeta(s+\alpha+\gamma)\zeta(s+\beta+\delta) \frac{\left( \frac {4\pi^2XMNv_1v_2}{T^2}\right)^s}{s}~ds
$$
and
$$ \int_0^\infty v_1^{s+\alpha+\gamma-1}(e(tv_1)+e(-tv_1))~dv_1 =t^{-s-\alpha-\gamma}\chi(1-s-\alpha-\gamma),$$ 
and similarly for the integral over $v_2$.
Incorporating these, we have simplified things to  
\begin{eqnarray*}&& \big(\frac {T}{2\pi}\big)^{-\alpha-\beta-\gamma-\delta}
\sum_{M\le N\le Q\atop (M,N)=1}\phi(\frac MQ)\phi(\frac NQ) M^{\gamma+\beta-1}N^{\delta+\alpha-1}\\
&&\qquad 
 \times \int_0^\infty \psi(t) t^{-\alpha-\beta-\gamma-\delta}
\frac{1}{2\pi i} \int_{(2)}\zeta(1-s-\alpha-\gamma)\zeta(1-s-\beta-\delta) \frac{ ( \frac {4\pi^2XMN }{t^2T^2} )^s}{s}~ds
~dt .
\end{eqnarray*}

The above expression is unchanged if $(\alpha,\gamma)$ is interchanged with  $(\beta,\delta)$. So the result of summing
terms for which $n_1/m_1\le 1$ rather than $m_1/n_1\le 1$ allows for 
 summing  over $M$ and $N$ as in Section 9; we obtain
\begin{eqnarray*}&& \big(\frac{T}{2\pi}\big)^{-\alpha-\beta-\gamma-\delta}\frac{1}{(2\pi i)^2} \int_{z,w} \tilde{\phi}(z)\tilde{\phi}(w)
\frac{\zeta(1-\beta-\gamma-s+z)\zeta(1-\alpha-\delta-s+w)}{\zeta(2-\alpha-\beta-\gamma-\delta-2s+z+w)}
  \\
&&\qquad 
 \times \int_0^\infty \psi(t) t^{-\alpha-\beta-\gamma-\delta}
\frac{1}{2\pi i} \int_{(2)}\zeta(1-s-\alpha-\gamma)\zeta(1-s-\beta-\delta) \frac{ ( \frac {X }{t^2T^2} )^sQ^{z+w}}{s}~ds~dw~dz
~dt .
\end{eqnarray*}
Moving the $s$-path to the right to $\infty$ and the $z$ and $w$ paths to the left to $-1/4$, say we obtain
\begin{eqnarray*}&&
\int_0^\infty \psi(t) \bigg(\left(\frac{tT}{2\pi}\right)^{-\alpha-\beta-\gamma-\delta}\frac{\zeta(1-\alpha-\gamma)\zeta(1-\beta-\delta)\zeta(1-\beta-\gamma)
\zeta(1-\alpha-\delta)}{\zeta(2-\alpha-\beta-\gamma-\delta)}\\&&
\qquad  +  X^{-\alpha-\gamma}\left(\frac{tT}{2\pi}\right)^{\alpha-\beta+\gamma-\delta} \frac{\zeta(1+\alpha-\beta+\gamma-\delta)\zeta(1+\alpha-\beta)
\zeta(1+\gamma-\delta)}{(\alpha+\gamma)~\zeta(2+\alpha-\beta+\gamma-\delta)}\\
&&\qquad \qquad + X^{-\beta-\delta}\left(\frac{tT}{2\pi}\right)^{-\alpha+\beta-\gamma+\delta} 
\frac{\zeta(1-\alpha+\beta-\gamma+\delta)\zeta(1-\alpha+\beta)
\zeta(1-\gamma+\delta)}{(\beta+\delta)~ \zeta(2-\alpha+\beta-\gamma+\delta) }\\
&&\qquad + X^{-\alpha-\delta} \left(\frac{tT}{2\pi}\right)^{\alpha-\beta -\gamma+\delta}
\frac{ \zeta(1-\gamma+\delta)\zeta(1+\alpha-\beta-\gamma+\delta)\zeta(1+\alpha-\beta )}
{(\alpha+\delta)\zeta(2+\alpha-\beta-\gamma+\delta)}\\
&&\qquad +  X^{-\beta-\gamma} \left(\frac{tT}{2\pi}\right)^{-\alpha+\beta +\gamma-\delta}
\frac{ \zeta(1-\alpha+\beta)\zeta(1-\alpha+\beta+\gamma-\delta)\zeta(1+\gamma-\delta)}
{(\beta+\gamma)\zeta(2-\alpha+\beta+\gamma-\delta)}
\bigg)~dt
\end{eqnarray*}
with an error term of $O(Q^{-1/4})$.
This expression is exactly what we were hoping for; it is
identical to the ``two-swap'' terms found in the descent approach of section 9 of [CK].

\section{Conclusion}
We have shown how to reproduce the complete conjecture for the shifted fourth moment of $\zeta$ by analyzing the mean square of 
long Dirichlet polynomials whose coefficients are convolutions of two smooth arithmetic functions. In the next paper we will carry this analysis 
out for coefficients which are convolutions of an arbitrary number of convolutions and use this to reproduce the full conjecture 
for the $2k$th moment of $\zeta$ for an arbitrary $k$.

\end{document}